\documentclass[reqno,centertags, 12pt]{amsart}
\usepackage{amsmath,amsthm,amscd,amssymb,amsfonts}
\usepackage{latexsym,verbatim}

-\marginparwidth \oddsidemargin -4mm \evensidemargin \oddsidemargin

\newcommand{\bbR}{{\mathbb{R}}}

\newcommand{\bbT}{{\mathbb{T}}}

\newcommand{\calI}{{\mathcal I}}

\newcommand{\lb}{\label}

\newcommand{\beq}{\begin{equation}}
\newcommand{\eeq}{\end{equation}}
\newcommand{\ba}{\begin{align}}
\newcommand{\ea}{\end{align}}
\newcommand{\eps}{\varepsilon}

\newcommand{\til}{\tilde}

\newcommand{\ffi}{\varphi}

\newcommand{\cea}{c_e^*(A)}
\newcommand{\dea}{D_e(A)}

\newcounter{smalllist}

 \DeclareMathOperator{\Ima}{Im}

 \DeclareMathOperator*{\sgn}{sgn}

\allowdisplaybreaks \numberwithin{equation}{section}

\newtheorem{theorem}{Theorem}[section]

\newtheorem{lemma}[theorem]{Lemma}

\theoremstyle{definition}

\theoremstyle{remark}

\begin{document}
\title[Pulsating Front Speed-up and Diffusion Enhancement]
{Sharp Asymptotics for KPP Pulsating Front Speed-up and Diffusion
Enhancement by Flows}

\author{Andrej Zlato\v s}


\thanks{Department of Mathematics, University of Chicago,
Chicago, IL 60637; email: {zlatos@math.uchicago.edu}}
\thanks{The author acknowledges partial support by the NSF through the grant DMS-0632442}

\begin{abstract}
We study KPP pulsating front speed-up and effective diffusivity
enhancement by general periodic incompressible flows. We prove the
existence of and determine the limits $c^*(A)/A$ and $D(A)/A^2$ as
$A\to\infty$, where $c^*(A)$ is the minimal front speed and $D(A)$
the effective diffusivity.
\end{abstract}

\maketitle

\section{Introduction} \lb{S1}

We study reaction-diffusion fronts in the presence of strong
incompressible flows. We consider the PDE
\begin{equation}\label{1.1}
T_t+Au\cdot\nabla T=\Delta T+f(T)
\end{equation}
on $\bbR^n$, with $T(t,x)\in[0,1]$ the normalized temperature of a
premixed combustible gas. The non-linear reaction rate $f$ is of
Kolmogorov-Petrovskii-Piskunov (KPP) type \cite{KPP}:
\begin{gather}
 f\in C^{1,\eps}([0,1]), \notag \\
 \text{$f(0)=f(1)=0$ and $f$ is non-increasing on $(1-\eps,1)$ for some $\eps>0$}, \label{1.0b} \\
 0<f(s)\le sf'(0) \text{ for $s\in(0,1)$}. \notag
\end{gather}
The 1-periodic flow $u:\bbT^n\to\bbR^n$ satisfies
\begin{equation}\label{1.0a}
 u\in C^{1,\eps}(\bbT^n), \qquad
 \nabla\cdot u \equiv 0, \qquad
 \int_{\bbT^n} u\,dx = 0.
\end{equation}
That is, $u$ is incompressible and mean-zero.

The number $A\in\bbR$ is the flow amplitude. We will consider the
case of strong flows (i.e., large $A$) and their influence on the
speed of propagation of pulsating fronts for \eqref{1.1}. This
problem has recently seen increased activity and has been addressed
by various authors
--- see, e.g., \cite{ABP,B,BHN-2,CKOR,Heinze,KR,NR,RZ}.

A {\it pulsating front} in the direction $e\in\bbR^n$, $|e|=1$, is a
solution of \eqref{1.1} of the form $T(t,x)=U(x\cdot e-ct,x)$, with
$c$ the front speed, and $U$ 1-periodic in $x$ and such that
\begin{eqnarray*} 
&&\lim_{s\to-\infty}U(s,x)=1,\\
&&\lim_{s\to+\infty}U(s,x)=0,
\end{eqnarray*}
uniformly in $x$. It is well known \cite{BHN-1} that in the KPP case
there is $\cea$, called the {\it minimal pulsating front speed},
such that pulsating fronts exist precisely for $c\ge \cea$ (we
suppress the $u$ and $f$ dependence in our notation). We note that
$\cea$ also determines the propagation speed of solutions to the
Cauchy problem with general compactly supported initial data
\cite{BHN-1, Weinberger}.

Mixing by flows (coupled to diffusion) typically increases the speed
of pulsating fronts for \eqref{1.1}. The minimal front speed $\cea$
can grow at most linearly with $A$ \cite{BHN-2} and does so for
shear (unidirectional) flows \cite{ABP,B,CKOR,Heinze}
\begin{equation}\label{1.1a}
u(x)=(\alpha(x'),0,\dots,0) \quad(x'=(x_2,\dots,x_n))
\qquad\text{and}\qquad e=(1,0,\dots,0).
\end{equation}
The same is true for so-called {\it percolating} flows which possess
infinite channels \cite{CKOR}, contrasting with the case of {\it
cellular} flows when, at least in two dimensions, $\cea=O(A^{1/4})$
\cite{ABP,CKOR,Heinze,NR} (see also \cite{RZ} for a
three-dimensional example).

We are interested here in all flows which maximally (i.e., linearly)
enhance the minimal front speed for \eqref{1.1} and our goal is to
determine the asymptotic rate of this front speed-up --- to prove
the existence and evaluate the limit of $\cea/A$ as $A\to\infty$.
For shear flows, this limit has been known to exist \cite{B} and has
been determined in \cite{Heinze}, but both problems have been open
in general.

We thus consider general periodic flows \eqref{1.0a} and let
\begin{equation}\label{1.2}
\calI \equiv \big\{ w\in H^1(\bbT^n) \,\,\,\big|\,\,\, \Ima w=0
\text{ and } u\cdot\nabla w = 0 \big\}
\end{equation}
be the set of {\it real-valued first integrals} of the flow $u$. We
then have the following main result.


\begin{theorem}\lb{T.1.1}
If  $u$ and $f$ satisfy \eqref{1.0b} and \eqref{1.0a} and
$|e|=1$, then
\begin{equation}\label{1.4a}
\lim_{A\to\infty} \frac \cea A = \sup_{\substack{w\in\calI \\
\|\nabla w\|_2^2\le f'(0)\|w\|_2^2}} \frac{\int_{\bbT^n}(u\cdot
e)w^2\,dx}{\|w\|_2^2}.
\end{equation}
In particular, the limit exists.  Moreover,
\begin{align}
& \lim_{f'(0)\to 0} \lim_{A\to\infty} \frac \cea {2\sqrt{f'(0)}A} =
\sup_{w\in\calI} \frac{\int_{\bbT^n}(u\cdot e)w\,dx}{\|\nabla
w\|_2}, \label{1.6} \\
& \lim_{f'(0)\to \infty} \lim_{A\to\infty} \frac \cea {A}  =
\sup_{w\in\calI} \frac{\int_{\bbT^n}(u\cdot e)w^2\,dx}{\|w\|_2^2}
\le \max_{x\in\bbT^n} \{u(x)\cdot e\}. \label{1.6a}
\end{align}
\end{theorem}

{\it Remarks.} 1. Inequality ``$\ge$'' in  \eqref{1.4a} (with
$\liminf_{A\to\infty}$ in place of $\lim_{A\to\infty}$) has been
proved in \cite{BHN-2}, and \cite{Heinze} showed equality in the
case of shear flows \eqref{1.1a}.
\smallskip

2. \eqref{1.6a} already appeared in \cite{BHN-2}, with either
$\liminf_{A\to\infty}$ or $\limsup_{A\to\infty}$ in place of
$\lim_{A\to\infty}$. For shear flows \eqref{1.1a} the inequality
becomes an equality \cite{Heinze} due to \eqref{1.4a} and continuity
of $u$.
\smallskip

3. Notice that \eqref{1.4a} (for any $f'(0)$) is positive precisely
when there exists $w\in\calI$ such that $\int_{\bbT^n}(u\cdot
e)w\,dx\neq 0$ (take $1\pm \eps w$ in \eqref{1.4a}). This is also
the condition for  positivity of \eqref{1.6} and \eqref{1.12} below.
\smallskip

4. The result extends directly to the more general case of
$x$-dependent and 1-periodic reaction and second-order term (see
Theorem \ref{T.3.3}). We perform the proof in the simpler setting
above for the sake of transparency.
\smallskip

It has been shown in \cite{NR,RZ} that, at least in two dimensions,
there is a close relationship between the minimal front speeds for
\eqref{1.1} and the effective diffusivity in the homogenization
theory for the related advection-diffusion problem
\begin{equation}\label{1.7}
\Phi_t+Au\cdot\nabla \Phi=\Delta \Phi.
\end{equation}
As is well known, the long-time behavior of solutions to \eqref{1.7}
is governed by the effective diffusion equation
\begin{equation*}
\Psi_t=\sum_{i,j=1}^n \sigma_{ij}(A)\frac{\partial^2\Psi}{\partial
x_i\partial x_j}.
\end{equation*}
Here $\sigma(A)$ is a constant effective diffusivity matrix. If
$e\in\bbR^n$ and we let $\chi_{e,A}$ be the mean-zero solution of
\begin{equation}\label{1.9}
-\Delta\chi_{e,A}+Au\cdot\nabla\chi_{e,A}=Au\cdot e
\end{equation}
on $\bbT^n$, then $\sigma(A)$ is given by
\begin{equation*}
e\cdot\sigma(A)e' = \int_{\bbT^n}(\nabla\chi_{e,A} +
e)\cdot(\nabla\chi_{e',A} + e') dx = e\cdot e' +
\int_{\bbT^n}\nabla\chi_{e,A}\cdot\nabla\chi_{e',A} dx.
\end{equation*}
The {\it effective diffusivity} for \eqref{1.7} in the direction
$e\in\bbR^n$, $|e|=1$, is now
\begin{equation}\label{1.11}
\dea\equiv e\cdot\sigma(A)e = 1+\|\nabla\chi_{e,A}\|_2^2.
\end{equation}

Again, mixing by flows enhances the effective diffusivity. It is
easy to show that $\dea$ can grow at most quadratically with $A$,
and flows that achieve this are said to {\it maximally enhance
diffusion} (see \cite{BGW,FP,KM} and references therein). It turns
out that our method applies to the problem of determining the
asymptotic rate of this enhancement as well, and we find the limit
$\dea/A^2$ as $A\to\infty$ for general periodic flows. To the best
of the author's knowledge, existence of this limit has not been
known before.

\begin{theorem}\lb{T.1.2}
If $u$ satisfies \eqref{1.0a} and
$|e|=1$, then
\begin{equation}\label{1.12}
\lim_{A\to\infty} \frac \dea {A^2} = \sup_{w\in\calI}
\bigg(\frac{\int_{\bbT^n}(u\cdot e)w\,dx}{\|\nabla w\|_2}\bigg)^2.
\end{equation}
In particular, the limit exists. Moreover, there is $w_0\in\calI$
which is a maximizer of \eqref{1.12} and $\chi_{e,A}/A\to w_0$ in
$H^1(\bbT^n)$.
\end{theorem}

{\it Remarks.} 1. It follows that the left hand side of \eqref{1.12}
is the square of the left hand side of \eqref{1.6}. This has been
established in two dimensions by Ryzhik and the author \cite{RZ},
even without the $A\to\infty$ limit (see also \cite{NR}).
\smallskip

2. We show that  if  \eqref{1.12} is positive, then the maximizers
are precisely $w=aw_0+b$ with $a,b\in\bbR$, $a\neq 0$.
\smallskip

3. If one considers the small diffusion problem $\phi_t = \eps\Delta
\phi + u\cdot\nabla \phi$ instead of \eqref{1.7}, then the
corresponding effective diffusivity satisfies $\til D_e(\eps) = \eps
D_e(\eps^{-1})$. Hence the limit $\lim_{\eps\to 0} \eps \til
D_e(\eps)$ also equals \eqref{1.12}.
\smallskip

4. Again, there is a straightforward extension to the case of
$x$-dependent second order term and even non-mean-zero flows (see
Theorem \ref{T.2.1}).
\smallskip


We prove Theorem \ref{T.1.2} in Section \ref{S2} and Theorem
\ref{T.1.1} in Section \ref{S3}. The generalizations to the case of
$x$-dependent second-order and reaction terms are
Theorems~\ref{T.2.1} and \ref{T.3.3} below.

\section{Effective Diffusivity Enhancement} \lb{S2}

\begin{proof}[Proof of Theorem \ref{T.1.2}]
Let $\psi_A\equiv\chi_{e,A}/A$, so that
\begin{equation}\label{2.1}
-\Delta\psi_A+Au\cdot\nabla\psi_A=u\cdot e
\end{equation}
Multiplying this by $\psi_A$ and integrating over $\bbT^n$ we obtain
using incompressibility of the flow,
\begin{equation}\label{2.2}
\|\nabla\psi_A\|_2^2 = \int_{\bbT^n} (u\cdot e)\psi_A dx \le
\|u\cdot e\|_2 \|\psi_A\|_2.
\end{equation}
Poincar\' e inequality
\begin{equation}\label{3.11}
\|w\|_2\le C\|\nabla w\|_2
\end{equation}
for some $C<\infty$ and any mean-zero $w$ then yields
\begin{equation}\label{2.3}
\|\psi_A\|_{H^1}  \le C\|u\cdot e\|_2.
\end{equation}
It also follows from \eqref{2.2} that
\begin{equation}\label{2.4}
\frac{\dea}{A^2} = \frac 1{A^2}+\|\nabla\psi_A\|_2^2 = \frac 1{A^2}
+ \int_{\bbT^n} (u\cdot e)\psi_A dx.
\end{equation}

Since $\|\psi_A\|_{H^1}$ is uniformly bounded, there is a sequence
$A_k\to\infty$ such that $\psi_{A_k}$ converges to some $w_0\in
H^1(\bbT^n)$, weakly in $H^1(\bbT^n)$ and strongly in $L^2(\bbT^n)$.
Then $\Delta\psi_{A_k}\to \Delta w_0$ and $\nabla\psi_{A_k}\to
\nabla w_0$ in the sense of distributions and \eqref{2.1} divided by
$A_k$ implies
\begin{equation}\label{2.5}
u\cdot\nabla w_0=0
\end{equation}
in the sense of distributions. Since $w_0\in H^1(\bbT^n)$, this
equality holds almost everywhere and $w_0\in\calI$. 
We also have
\[
\|\nabla w_0\|_2^2 \le \limsup_{k\to\infty} \|\nabla
\psi_{A_k}\|_2^2 = \int_{\bbT^n} (u\cdot e)w_0 dx = \int_{\bbT^n}
\nabla\psi_{A}\nabla w_0 dx \le  \|\nabla \psi_{A}\|_2 \|\nabla
w_0\|_2
\]
where we used \eqref{2.2} in the second step, and \eqref{2.1}
multiplied by $w_0$ and integrated over $\bbT^n$ (together with
\eqref{2.5}) in the third step. Thus
\begin{equation}\label{2.6}
\|\nabla w_0\|_2\le \|\nabla \psi_{A}\|_2
\end{equation}
as well as
\[
\limsup_{k\to\infty}\|\nabla \psi_{A_k}\|_2\le \|\nabla w_0\|_2.
\]
These give
\begin{equation*}
\lim_{k\to\infty}\|\nabla \psi_{A_k}\|_2 = \|\nabla w_0\|_2,
\end{equation*}
which turns the weak $H^1$-convergence into a strong one:
\begin{equation}\label{2.8}
\psi_{A_k}\to  w_0 \qquad\text{in $H^1(\bbT^n)$}.
\end{equation}

Let us assume $w_0\not\equiv 0$. Then $\nabla w_0\not\equiv 0$
because each $\psi_A$ is mean-zero. From \eqref{2.4} and
\eqref{2.8},
\begin{equation}\label{2.9}
\lim_{k\to\infty} \frac {D_e(A_k)} {A_k^2} = \|\nabla w_0\|_2^2 =
\int_{\bbT^n}(u\cdot e)w_0 \,dx = \bigg( \frac {\int_{\bbT^n}(u\cdot
e)w_0 \,dx} {\|\nabla w_0\|_2} \bigg)^2.
\end{equation}
Pick an arbitrary non-constant $w\in \calI$. If we multiply
\eqref{2.1} by $w$ and integrate, we obtain
\begin{equation}\label{2.10}
\bigg| \int_{\bbT^n}(u\cdot e)w \,dx \bigg| = \bigg|
\lim_{k\to\infty} \int_{\bbT^n} \nabla\psi_{A_k} \nabla w \,dx
\bigg| = \bigg| \int_{\bbT^n} \nabla w_0 \nabla w \,dx \bigg| \le
\|\nabla w_0\|_2 \|\nabla w\|_2.
\end{equation}
Hence
\[
\bigg( \frac {\int_{\bbT^n}(u\cdot e)w\,dx} {\|\nabla w\|_2}
\bigg)^2 \le \|\nabla w_0\|_2^2 = \lim_{k\to\infty} \frac {D_e(A_k)}
{A_k^2},
\]
with equality precisely when $\nabla w$ is a multiple of $\nabla
w_0$ (and so $w=aw_0+b$). This also means that $w_0$ is a maximizer
for \eqref{1.12}.

If now $B_k\to\infty$ is any sequence, then as above we can find a
subsequence (which we again call $B_k$) such that $\psi_{B_k}\to
w_1\in\calI$. 
But then $w_1$ must also maximize \eqref{1.12}, thus $w_1=aw_0+b$.
Moreover, $b=0$ because $\psi_A$ are mean-zero, and \eqref{2.9} with
$B_k$ in place of $A_k$ forces $a=1$. Hence $\psi_A\to w_0$ in
$H^1(\bbT^n)$ and \eqref{1.12} follows.

Finally, if $w_0\equiv 0$ is the only limit point of $\psi_{A}$,
then $\psi_A\to 0$ in $H^1(\bbT^n)$, and \eqref{1.12} follows from
\eqref{2.4} and \eqref{2.10}.
\end{proof}

Notice that \eqref{2.4}, \eqref{2.6}, and \eqref{2.9} show that
$\dea\ge 1+\delta A^2$, where $\delta$ is the limit in \eqref{1.12}.

We also note that in the special case of shear flows
$u(x)=(\alpha(x'),0,\dots,0)$ equation \eqref{1.9} becomes
\[
-\Delta_{x'}\chi_{e,A}=Ae_1\alpha(x')
\]
with $\chi_{e,A}(x)=\chi_{e,A}(x')$. Hence
$\chi_{e,A}=Ae_1\nabla_{x'}(-\Delta_{x'})^{-1}\alpha$ and the limit
in \eqref{1.12} equals
$|e_1|\|\nabla_{x'}(-\Delta_{x'})^{-1}\alpha\|_2^2$. This can be
found, e.g.,  in \cite[Lemma 7.3]{FP}.

As mentioned above, the result easily extends to the case of
$x$-dependent second order term and a non-mean-zero flow. We
consider
\begin{equation}\label{2.11}
\Phi_t+Au\cdot\nabla \Phi=\nabla\cdot(a\nabla \Phi)
\end{equation}
instead of \eqref{1.7} with 1-periodic and real symmetric uniformly
elliptic matrix $a$ and 1-periodic flow $u$ such that
\begin{equation}\label{2.12}
a\in C^2(\bbT^n),\qquad u\in C^{1,\alpha}(\bbT^n), \qquad
 \nabla\cdot u \equiv 0, \qquad
 \bar u \equiv \int_{\bbT^n} u\,dx.
\end{equation}
Then \eqref{1.9} and \eqref{1.11} are replaced by
\begin{gather*}
 -\nabla\cdot(a\nabla\chi_{e,A})+Au\cdot\nabla\chi_{e,A}=A(u-\bar u)\cdot e, \\
 \dea \equiv |||\nabla\chi_{e,A}+e|||_2^2,
\end{gather*}
with $|||w|||_2^2\equiv\int_{\bbT^n}\nabla w\cdot(a\nabla w)\,dx$.
If we define
\[
\calI_0\equiv \bigg\{w\in\calI \,\bigg|\, \int_{\bbT^n} w\,dx=0
\bigg\},
\]
then we have

\begin{theorem}\lb{T.2.1}
If $a$ and $u$ satisfy \eqref{2.12} and
$|e|=1$, then
\begin{equation}\label{2.13}
\lim_{A\to\infty} \frac \dea {A^2} = \sup_{w\in\calI_0}
\bigg(\frac{\int_{\bbT^n}(u\cdot e)w\,dx}{|||\nabla w|||_2}\bigg)^2.
\end{equation}
In particular, the limit exists. Moreover, there is $w_0\in\calI_0$
which is a maximizer of \eqref{2.13} and $\chi_{e,A}/A\to w_0$ in
$H^1(\bbT^n)$.
\end{theorem}

\section{KPP Front Speed-up} \lb{S3}

In this section we prove Theorem \ref{T.1.1}. We start with an
auxiliary lemma. Let us define
\begin{equation}\label{1.3}
\kappa_e(\lambda) \equiv \sup_{w\in\calI} \bigg\{ \frac{\lambda
\int_{T^n} (u\cdot e)w^2 \,dx -\|\nabla w\|_2^2}{\|w\|_2^2} \bigg\}.
\end{equation}
Note that $\kappa_e(\lambda)$ must be convex as it is a supremum of
linear functions. Also, $\kappa_e(\lambda)\ge 0$ because $w\equiv
1\in\calI$.

\begin{lemma}\lb{L.3.1}
Assume the setting of Theorem \ref{T.1.1}. Then for each
$\lambda>0$, the supremum in \eqref{1.3} is attained, the maximizer
is unique up to multiplication, and
\begin{equation}\label{1.4}
\lim_{A\to\infty} \frac \cea A = \inf_{\lambda>0}
\frac{f'(0)+\kappa_e(\lambda)}{\lambda}.
\end{equation}
\end{lemma}

\begin{proof}
It has been shown in \cite{BHN-1} that the minimal front speed
$\cea$ can be computed using the variational principle
\begin{equation}\label{3.1}
\cea=\inf_{\lambda>0}\frac{f'(0)+\lambda^2+\kappa(\lambda;A)}{\lambda}.
\end{equation}
Here $\kappa(\lambda;A)$ is the unique eigenvalue of the problem
\begin{equation}\label{3.2}
\Delta\varphi-Au\cdot\nabla\varphi-2\lambda e\cdot\nabla\varphi+
{\lambda}Au\cdot e\varphi=\kappa(\lambda;A)\varphi,\qquad
\hbox{$\varphi>0$}
\end{equation}
on $\bbT^n$, with a unique normalized eigenfunction
$\varphi_A(x;\lambda)$. Moreover, the function
\[
\mu(\lambda;A) \equiv \lambda^2 + \kappa(\lambda;A)
\]
is monotonically increasing and convex in $\lambda\ge 0$, with
$\mu(0;A)=0$ (see \cite{BH,NR}).

We now rewrite \eqref{3.1} and \eqref{3.2} as
\begin{equation}\label{3.3}
\frac \cea A =
\inf_{\lambda>0}\frac{f'(0)+(\lambda/A)^2+\kappa(\lambda/A;A)}{\lambda}.
\end{equation}
and
\begin{equation}\label{3.4}
\Delta\varphi_A-Au\cdot\nabla\varphi_A- \frac{2\lambda} A
e\cdot\nabla\varphi_A+ {\lambda}u\cdot
e\varphi_A=\kappa(\lambda/A;A)\varphi_A,\qquad \hbox{$\varphi_A>0$}.
\end{equation}

We multiply \eqref{3.4} by $\ffi_A^{-1}$ and integrate to obtain
(using incompressibility of $u$)
\begin{equation}\label{3.5}
0\le  \|\nabla \ln\ffi_A\|_2^2 = \kappa(\lambda/A;A).
\end{equation}
Similarly, multiplication by $\ffi_A$ yields
\begin{equation}\label{3.6}
\kappa(\lambda/A;A) + \|\nabla \ffi_A\|_2^2 = \lambda\int_{\bbT^n}
(u\cdot e)\ffi_A^2 dx \le \lambda \|u\cdot e\|_\infty .
\end{equation}
since $\|\ffi_A\|_2=1$. This again means that there is a sequence
$A_k\to\infty$ such that $\ffi_{A_k}$ converges to some $w_0\in
H^1(\bbT^n)$, weakly in $H^1(\bbT^n)$ and strongly in $L^2(\bbT^n)$.
The convergence $\Delta\ffi_{A_k}\to\Delta w_0$ and
$\nabla\ffi_{A_k}\to\nabla w_0$ in the sense of distributions,
boundedness of $\kappa(\lambda/A;A)$ in $A$, and \eqref{3.4} divided
by $A$ then imply \eqref{2.5}
and so $w_0\in\calI$ (note that $\|w_0\|_2=\|\ffi_{A_k}\|_2=1$).

Now we multiply \eqref{3.4} by $w_0$ and integrate to obtain (with
$o(1)=o(k^0)$ and using \eqref{3.6})
\begin{align*}
-\int_{\bbT^n} \nabla \ffi_{A_k} \nabla w_0 dx + \lambda
\int_{\bbT^n} (u\cdot e) w_0^2 dx + o(1) & =
\kappa(\lambda/A_k;A_k) + o(1) \\
& = \lambda \int_{\bbT^n} (u\cdot e) w_0^2 dx
-\|\nabla\ffi_{A_k}\|_2^2 + o(1).
\end{align*}
Once again it follows that
\[
\|\nabla w_0\|_2^2 \le \limsup_{k\to\infty} \|\nabla
\ffi_{A_k}\|_2^2 \le \|\nabla w_0\|_2 \limsup_{k\to\infty}\|\nabla
\ffi_{A_k}\|_2
\]
and so as in Section \ref{S2},
\begin{equation}\label{3.8}
\ffi_{A_k}\to  w_0 \qquad\text{in $H^1(\bbT^n)$}.
\end{equation}
\eqref{3.6} then yields
\begin{equation*}
\kappa_0 \equiv \lim_{k\to\infty} \kappa(\lambda/A_k;A_k) =
\lambda\int_{\bbT^n} (u\cdot e)w_0^2 dx -\int_{\bbT^n} |\nabla
w_0|^2\, dx.
\end{equation*}

Let $w\in\calI\cap L^\infty(\bbT^n)$, multiply \eqref{3.4} for
$A=A_k$ by $w^2/\ffi_{A_k}$ and integrate to obtain (using that
$\nabla\ffi_A/\ffi_A=\nabla\ln\ffi_A$ are uniformly bounded in
$L^2(\bbT^n)$ by \eqref{3.5} and \eqref{3.6})
\begin{align}
\kappa_0 \|w\|_2^{2} & = \lambda\int_{\bbT^n} (u\cdot e)w^2 dx
+\lim_{k\to\infty} \int_{\bbT^n} \bigg|\frac{\nabla
\ffi_{A_k}}{\ffi_{A_k}}\bigg|^2 w^2 - 2 \frac{\nabla
\ffi_{A_k}}{\ffi_{A_k}} w\nabla w
- \frac{2\lambda}{A_k} e\cdot \frac{\nabla \ffi_{A_k}}{\ffi_{A_k}} w^2\,dx  \notag \\
& \ge \lambda\int_{\bbT^n} (u\cdot e)w^2 dx - \int_{\bbT^n} |\nabla
w|^2 \,dx. \label{3.9a}
\end{align}
Since each $w\in\calI$ is the $H^1$-limit of $w_N(x)\equiv
\sgn(w(x))\min\{|w(x)|,N\}\in\calI\cap L^\infty(\bbT^n)$, this
inequality extends to all $w\in\calI$. Hence
$\kappa_0=\kappa_e(\lambda)$ from \eqref{1.3}, and $w_0$ is a
maximizer for \eqref{1.3} (because $\|w_0\|_2=1$). Moreover, if
$B_k\to\infty$ is any sequence with
\[
\lim_{k\to\infty} \kappa(\lambda/B_k;B_k) \equiv \kappa_1 ,
\]
then repeating the above argument we find that there must be a
subsequence (which we again call $B_k$) such that $\ffi_{B_k}\to
w_1\in\calI$ in $H^1(\bbT^n)$, $\|w_1\|_2=1$. But then as before,
\[
\lambda\int_{\bbT^n} (u\cdot e)w_1^2 dx - \int_{\bbT^n} |\nabla
w_1|^2 \,dx = \kappa_1 \ge \frac{\lambda\int_{\bbT^n} (u\cdot e)w^2
dx - \int_{\bbT^n} |\nabla w|^2 \,dx}{\|w\|_2^2}
\]
for any $w\in\calI$. Taking $w=w_0$ we obtain
$\kappa_1=\kappa_e(\lambda)$, and so
\begin{equation}\label{3.10}
\kappa_e(\lambda) = \lim_{A\to\infty} \kappa(\lambda/A;A) =
\lambda\int_{\bbT^n} (u\cdot e)w_0^2 dx -\int_{\bbT^n} |\nabla
w_0|^2\, dx.
\end{equation}

The function $\kappa_e(\lambda)$ is  convex, monotonically
increasing, and non-negative, as it is the pointwise limit of
functions $\mu(\lambda/A;A)=(\lambda/A)^2 + \kappa(\lambda/A;A)$
which have the same properties. This also implies that the
convergence in \eqref{3.10} is uniform on each bounded interval of
$\lambda$. We then have
\[
\lim_{A\to\infty}
\inf_{\lambda>0}\frac{f'(0)+\mu(\lambda/A;A)}{\lambda}=
\inf_{\lambda>0} \frac{f'(0)+\kappa_e(\lambda)}{\lambda}
\]
($\le$ is immediate, whereas $\ge$ uses convexity of
$\mu(\lambda/A;A)$ once more). This proves \eqref{1.4}.

We are left with showing that any maximizer of \eqref{1.3} is a
multiple of $w_0$. Denote $\ffi_k\equiv\ffi_{A_k}$ and notice that
\eqref{3.8} shows that (after passing to a subsequence --- we will
repeat this without mentioning it below), $\nabla\ffi_k(x)\to\nabla
w_0(x)$ and $\ffi_k(x)\to w_0(x)$ for a.e. $x$. Next \eqref{3.5} and
\eqref{3.11} imply that if $c_k$ is the average of $\ln\ffi_k$, then
$\ln\ffi_k-c_k\to\omega$ strongly in $L^2$ and weakly in $H^1$. But
then $\ln \ffi_k(x)-c_k\to\omega(x)$ for a.e. $x$. Since $\ln
\ffi_k(x)\to \ln w_0(x)$ for a.e. $x$, it follows that $c_k\to c$
and $\omega=\ln w_0-c$. We thus obtain $\ln w_0\in H^1$ which means
$w_0(x)\neq 0$ for a.e. $x$, and so for a.e. $x$,
\begin{equation}\label{3.10a}
\frac{\nabla \ffi_k(x)}{\ffi_k(x)}\to\frac{\nabla w_0(x)}{w_0(x)}.
\end{equation}

Let now $w\not\equiv 0$ be a maximizer of \eqref{1.3} and let us
first assume $w\ge 0$ almost everywhere. Then \eqref{3.9a} for $w_N$
and $w_N\to w$ in $H^1$ show
\[
\lim_{N\to\infty} \limsup_{k\to\infty} \bigg\|\frac{\nabla
\ffi_k}{\ffi_k}w_N-\nabla w_N \bigg\|_2 =0.
\]
But then \eqref{3.10a} and pointwise convergence of $w_N$ and
$\nabla w_N$ to $w$ and $\nabla w$, respectively, give for a.e. $x$,
\[
\frac{\nabla w_0(x)}{w_0(x)}w(x)=\nabla w(x).
\]
We now let $w_\eps(x)\equiv \max\{w(x),\eps\}$ so that $\ln
w_\eps\in H^1$ and
\[
\nabla\ln w_\eps(x) = \begin{cases} \nabla \ln w_0(x) &
w_\eps(x)>\eps, \\ 0 & w_\eps(x)=\eps.
\end{cases}
\]
This and $\ln w_0\in H^1$ means that $\|\nabla \ln w_\eps\|_2$ is
bounded, and again we must have $\ln w_{\eps_k}-c_{\eps_k}\to\omega$
strongly in $L^2$, weakly in $H^1$, and pointwise almost everywhere.
But $\ln w_{\eps_k} (x)\to \ln w(x)$, so again $c_{\eps_k}\to c$ and
$\ln w\in H^1$. Hence $w(x)> 0$ for a.e. $x$, and so $\nabla\ln w(x)
= \nabla \ln w_0(x)$ for a.e. $x$. This means $\ln w-\ln w_0$ is
constant, that is, $w$ is a multiple of $w_0$.

If $w$ is an arbitrary maximizer of \eqref{1.3}, then both
$w_\pm(x)\equiv\max\{\pm w(x),0\}\in\calI$ must be maximizers of
\eqref{1.3} (or $\equiv 0$). But then $w_\pm(x)>0$ for a.e. $x$,
meaning that one of them is zero while the other is a multiple of
$w_0$.
\end{proof}

\begin{proof}[Proof of Theorem \ref{T.1.1}]
Inequality ``$\ge$'' in \eqref{1.4a} is immediate from \eqref{1.3}
and \eqref{1.4}. To prove the opposite inequality it is sufficient
to find $\lambda$ such that the unique normalized non-negative
maximizer $w_0(\lambda)$ of \eqref{1.3} satisfies
$\gamma(\lambda)\equiv \|\nabla w_0(\lambda)\|_2^2=f'(0)$.

To this end notice that if $\lambda=0$, then $w_0(\lambda)\equiv 1$
and so $\gamma(0)=0$. Also, $\gamma$ must be continuous. Indeed ---
let $\lambda_k\to\lambda_\infty<\infty$ and denote $w_k\equiv
w_0(\lambda_k)$. Then \eqref{3.6} and \eqref{3.8} imply that $w_k$
are uniformly bounded in $H^1$. Thus a subsequence (again called
$w_k$) converges strongly in $L^2$ and weakly in $H^1$ to some
$\omega$. Obviously $\omega\in\calI$ as well as
$\lambda_k\int(u\cdot e)w_k^2\,dx\to \lambda_\infty\int(u\cdot
e)\omega^2\,dx$ and $\|\nabla\omega\|_2\le\liminf \|\nabla w_k\|_2$.
But $\|\nabla\omega\|_2<\liminf \|\nabla w_k\|_2$ is impossible
(otherwise $w_k$ would not maximize \eqref{1.3} for large $k$) and
so for a subsequence,
\[
\kappa_e(\lambda_k) = \lambda_k\int(u\cdot e)w_k^2\,dx - \|\nabla
w_k\|_2^2 \to \lambda_\infty\int(u\cdot e)\omega^2\,dx -
\|\nabla\omega\|_2^2.
\]
Since $\kappa_e$ is continuous, this means $\omega=w_\infty$. We
have thus proved that every sequence $\lambda_k\to\lambda_\infty$
has a subsequence with
$\gamma(\lambda_{k_j})\to\gamma(\lambda_\infty)$, that is, $\gamma$
is continuous with $\gamma(0)=0$.

Let now $\Gamma\equiv\sup_{\lambda\ge 0}\gamma(\lambda)$. If
$f'(0)\in(0,\Gamma)$, then there is $\lambda>0$ with
$\gamma(\lambda)=f'(0)$ and \eqref{1.4a} is proved. If, on the other
hand, $\Gamma<\infty$ and $f'(0)\ge\Gamma$, then \eqref{1.4} is
bounded above by
\[
\liminf_{\lambda\to\infty} \frac {\lambda\int(u\cdot
e)w_0(\lambda)^2\,dx - \|\nabla w_0(\lambda)\|_2^2 + f'(0)}
{\lambda} = \liminf_{\lambda\to\infty} \int(u\cdot
e)w_0(\lambda)^2\,dx,
\]
which does not exceed the right hand side of \eqref{1.4a} due to
$\|\nabla w_0(\lambda)\|_2^2\le \Gamma\le f'(0)
\|w_0(\lambda)\|_2^2$.


Since \eqref{1.6a} is immediate from \eqref{1.4a}, we are left with
proving \eqref{1.6}. Let us consider any $w\in\calI$ with
$\|w\|_2=1$ and $\|\nabla w\|_2^2\le f'(0)$. Let $\bar w\equiv
\int_{\bbT^n} w\,dx\in [-1,1]$ and $\omega\equiv w-\bar w$. Then
$\|\omega\|_2^2 \le C\|\nabla\omega\|_2^2\le Cf'(0)$, and so
\[
1=\int_{\bbT^n} \bar w^2+2\bar w\omega +\omega^2\,dx = \bar w^2 +
O(\sqrt{f'(0)})
\]
as $f'(0)\to 0$. Hence $\bar w=1+O(\sqrt{f'(0)})$ and we have
\[
\int_{\bbT^n}(u\cdot e)w^2\,dx = 2\int_{\bbT^n}(u\cdot e)\omega\,dx
+ O(f'(0)) \le 2\sqrt{f'(0)} \frac{\int_{\bbT^n}(u\cdot
e)\omega\,dx}{\|\nabla \omega\|_2} + O(f'(0))
\]
with equality when $\|\nabla \omega\|_2^2= f'(0)$. Picking first $w$
that maximizes \eqref{1.4a} and then $\omega$ that maximizes
\eqref{1.6} with $\|\nabla \omega\|_2^2= f'(0)$ (and adjusting $\bar
w$ accordingly) finishes the proof.
\end{proof}

Note that in the case of shear flows \eqref{1.1a} equation
\eqref{3.4} becomes
\begin{equation}\label{3.13a}
\Delta_{x'}\varphi +
{\lambda}\alpha\varphi=\kappa(\lambda/A;A)\varphi,\qquad
\hbox{$\varphi>0$}.
\end{equation}
with $\varphi(x)=\varphi(x')$. As a result
$\kappa(\lambda/A;A)=\kappa(\lambda;1)=\kappa_e(\lambda)$ and
$w_0(\lambda)=\varphi$, and \eqref{3.3} shows that $\cea/A$ is
non-increasing. This has been proved in \cite{B}. If the limit is
$\gamma$, then \eqref{3.3} gives
\[
\frac\cea A - \gamma \le \frac{2\sqrt{f'(0)}}A
\]
(which has been already observed in \cite{Heinze}). Here one uses
convexity of $\kappa_e$ and $\kappa_e(0)=0$ to show that the infimum
in \eqref{3.3} is achieved at some $\lambda\le \lambda_A\equiv
\sqrt{f'(0)}A$, as well as
\[
\inf_{\lambda>0} \frac{f'(0)+\kappa_e(\lambda)}\lambda \ge \min
\bigg\{ \inf_{\lambda\in(0,\lambda_A)}
\frac{f'(0)+\kappa_e(\lambda)}\lambda,
\frac{\kappa_e(\lambda_A)}{\lambda_A} \bigg\}.
\]
Moreover, if the infimum in \eqref{1.4} is achieved at a finite
$\lambda$, then \eqref{3.3} gives that
\[
\frac\cea A - \gamma =O(A^2).
\]
This condition is satisfied for all $f'(0)<\Gamma$, where $\Gamma$
is from the proof of Lemma \ref{L.3.1}, that is, it is the supremum
over $\lambda>0$ of the $\dot{H}^1$ norms of the principal
eigenfunctions of \eqref{3.13a}. This is because of \eqref{1.4}, the
definition of $\kappa_e$, and the fact that
\[
\lim_{\lambda\to\infty} \frac{\kappa_e(\lambda)}\lambda =
\sup_{w\in\calI} \frac{\int_{\bbT^n}(u\cdot e)w^2\,dx}{\|w\|_2^2}.
\]

Finally, we note that $\Gamma<\infty$ is possible --- in the shear
flow case it holds when there is an open set $U\subseteq \bbT^{n-1}$
such that $\alpha(x')=\max_{\bbT^{n-1}} \alpha$ for all $x'\in U$.
Then any $w\in H^1(\bbT^n)$ supported on $\bbT\times U$ and
independent of $x_1$ belongs to $\calI$ and maximizes \eqref{1.4a}
whenever $f'(0)\ge\|\nabla w\|_2^2/\|w\|_2^2$. Thus the limit in
\eqref{1.4a} need not be strictly increasing with $f'(0)$ (which
happens precisely when $\Gamma<\infty$).

In the more general case when the second order term and the
non-linearity depend on $x$, we consider
\begin{equation}\label{3.14}
T_t+Au\cdot\nabla T=\nabla\cdot(a\nabla T)+f(x,T)
\end{equation}
with $a$ 1-periodic real symmetric uniformly elliptic matrix and $u$
1-periodic such that
\begin{equation}\label{3.15}
a\in C^2(\bbT^n),\qquad u\in C^{1,\eps}(\bbT^n), \qquad
 \nabla\cdot u \equiv 0, \qquad
 \int_{\bbT^n} u\,dx=0.
\end{equation}
The non-linearity $f$ is 1-periodic in $x$ and satisfies for some
$\eps>0$
\begin{gather}
 f\in C^{1,\delta}(\bbT^n\times[0,1]), \notag \\
 \text{$f(x,0)=f(x,1)=0$ and $f(x,\cdot)$ is non-increasing on $(1-\eps,1)$ for each $x\in\bbT^n$}, \label{3.16} \\
 0<f(x,s)\le sf'_s(x,0) \text{ for $(s,x)\in(0,1)\times\bbT^n$}. \notag
\end{gather}
We let $\zeta(x)\equiv f'_s(x,0)>0$ and
$\zeta_0\equiv\int_{\bbT^n}\zeta(x)\,dx$. Equations \eqref{3.1} and
\eqref{3.2} are then replaced by (see \cite{BHN-1})
\begin{gather*}
 \cea=\inf_{\lambda>0}\frac{\kappa(\lambda,f;A)}{\lambda}, \\
 \nabla\cdot(a\nabla \varphi)  -Au\cdot\nabla\varphi-2\lambda e\cdot
a\nabla\varphi + [{\lambda}Au\cdot e + \zeta + \lambda^2 e\cdot ae
-\lambda\nabla\cdot(ae)] \varphi=\kappa(\lambda,f;A)\varphi.
\end{gather*}

If we now define
\begin{equation*} 
\kappa_e(\lambda,f)  \equiv \sup_{w\in\calI} \bigg\{ \frac{
\int_{T^n} (\lambda u\cdot e+ \zeta)w^2 \,dx -|||\nabla
w|||_2^2}{\|w\|_2^2} \bigg\},
\end{equation*}
then \eqref{1.4} becomes
\begin{equation*}
\lim_{A\to\infty} \frac {c_e^*(A,f)} A = \inf_{\lambda>0}
\frac{\kappa_e(\lambda,f)}{\lambda}
\end{equation*}
and mimicking the above proofs one obtains the following extension
of Theorem \ref{T.1.1}.

\begin{theorem}\lb{T.3.3}
If $a$, $u$, and $f$ satisfy \eqref{3.15} and \eqref{3.16} and
$|e|=1$, then
\begin{equation*}
\lim_{A\to\infty} \frac \cea A = \sup_{\substack{w\in\calI \\
|||\nabla w|||_2^2\le \int_{\bbT^n}\zeta w^2\,dx}}
\frac{\int_{\bbT^n}(u\cdot e)w^2\,dx}{\|w\|_2^2}.
\end{equation*}
In particular, the limit exists.
Moreover,
\begin{align*}
\lim_{\alpha\to 0} \lim_{A\to\infty} \frac {c_e^*(A,\alpha f)}
{2\sqrt{\alpha\zeta_0}A} & = \sup_{w\in\calI}
\frac{\int_{\bbT^n}(u\cdot e)w\,dx}{|||\nabla w|||_2}, \\
\lim_{\alpha\to \infty} \lim_{A\to\infty} \frac {c_e^*(A,\alpha f)}
{A} & =  \sup_{w\in\calI} \frac{\int_{\bbT^n}(u\cdot
e)w^2\,dx}{\|w\|_2^2} \le \max_{x\in\bbT^n} \{u(x)\cdot e\}.
\end{align*}
\end{theorem}

\end{document}